\documentclass[final]{colt2023}

\usepackage{array} %
\usepackage{caption} %
\usepackage{float} %
\usepackage{tablefootnote} %
\usepackage{multirow} %
\usepackage{enumitem}
\usepackage{ifthen}

\usepackage{url} %
\usepackage{mathtools}
\usepackage{mathrsfs} 
\makeatletter\def\@seccntformat#1{\protect\makebox[0pt][r]{\csname the#1\endcsname\hspace{12pt}}}\makeatother
\definecolor{mycolor1}{rgb}{0.105882,0.619608,0.466667}
\definecolor{mycolor2}{rgb}{0.85098,0.372549,0.00784314}
\definecolor{mycolor3}{rgb}{0.458824,0.439216,0.701961}
\definecolor{mycolor4}{rgb}{0.905882,0.160784,0.541176}
\definecolor{mycolor5}{rgb}{0.4,0.65098,0.117647}
\definecolor{mycolor6}{rgb}{0.65098,0.462745,0.113725}
\definecolor{mycolor7}{rgb}{0.901961,0.670588,0.00784314}
\definecolor{mycolor8}{rgb}{0.4,0.4,0.4}
\definecolor{mycolor9}{rgb}{0.301961,0,0.294118}
\definecolor{mycolor10}{rgb}{0.0313725,0.25098,0.505882}

\DeclarePairedDelimiter{\ceil}{\lceil}{\rceil}

\usepackage{bigfoot}

\newcommand{\cutchunk}[1]{}

\newcommand{\removesafe}[1]{}

\newcommand{\calA}{\mathcal{A}}

\newcommand{\calE}{\mathcal{E}}

\newcommand{\calM}{\mathcal{M}}

\newcommand{\dist}{\mathrm{dist}}

\newcommand{\reals}{{\mathbb{R}}}

\newcommand{\T}{\mathrm{T}}

\newcommand{\xorigin}{x_{\mathrm{ref}}}

\newcommand{\aref}[1]{\hyperref[#1]{A\ref{#1}}}
\newtheorem{question}[theorem]{Open Question}

\title[Polynomial linearly-convergent method for g-convex optimization?]{Open Problem: Polynomial linearly-convergent method for geodesically convex optimization?}
\usepackage{times}

 \coltauthor{\Name{Christopher Criscitiello} \Email{christopher.criscitiello@epfl.ch}\AND
\Name{Nicolas Boumal} \Email{nicolas.boumal@epfl.ch}\\
\addr }

\coltauthor{%
 \Name{Christopher Criscitiello} \Email{christopher.criscitiello@epfl.ch}\\
 \addr Ecole Polytechnique F\'ed\'erale de Lausanne (EPFL), Institute of Mathematics
 \AND
 \Name{David Mart\'inez-Rubio} \Email{martinez-rubio@zib.de}\\
 \addr Zuse Institute Berlin and Technische Universit\"{a}t
 \AND
 \Name{Nicolas Boumal} \Email{nicolas.boumal@epfl.ch}\\
 \addr Ecole Polytechnique F\'ed\'erale de Lausanne (EPFL), Institute of Mathematics
}

\begin{document}

\maketitle

\begin{abstract}%
Let $f \colon \calM \to \reals$ be a Lipschitz and geodesically convex function defined on a $d$-dimensional Riemannian manifold $\calM$.
Does there exist a first-order deterministic algorithm which (a) uses at most $O(\mathrm{poly}(d) \log(\epsilon^{-1}))$ subgradient queries to find a point with target accuracy $\epsilon$, and (b) requires only $O(\mathrm{poly}(d))$ arithmetic operations per query?
In convex optimization, the classical ellipsoid method achieves this.
After detailing related work, we provide an ellipsoid-like algorithm with query complexity $O(d^2 \log^2(\epsilon^{-1}))$ and per-query complexity $O(d^2)$
for the limited case where $\calM$ has constant curvature (hemisphere or hyperbolic space).
We then detail possible approaches and corresponding obstacles for designing an ellipsoid-like method for general Riemannian manifolds.
\end{abstract}


\section{Introduction and background} 
The ellipsoid method (developed by~\citet{shor1977cut} and~\citet{yudin1976evaluation,yudin1976informational}) is a fundamental algorithm in convex optimization and computational complexity.\footnote{For a description and analysis of the ellipsoid method, see~\citep[Sec.~3.2]{nemirovskibook} or~\citep[Sec.~3.2]{nesterov2004introductory}.}
The first proof that linear programs can be solved in polynomial time used the ellipsoid method~\citep{Kha79}.

Let us briefly review how the method works.
Suppose we seek to minimize a Lipschitz convex function $f \colon \reals^d \to \reals$ constrained to a closed ball $\bar B(\xorigin, r)$ of radius $r$ and center $\xorigin$. 
At each step of the ellipsoid method, there is an ellipsoid $E_k \subset \reals^d$ which by construction contains the minimizers of $\min_{x \in \bar B(\xorigin, r)} f(x)$.  The ellipsoid method queries the center of the ellipsoid and receives a subgradient, which in turn determines a halfspace containing the minimizers.  The next ellipsoid $E_{k+1}$ is the minimum-volume ellipsoid containing the intersection of the halfspace and $E_k$.  
The volume of the ellipsoids $E_k$ decreases at a linear rate, leading to the linear rate of convergence for the method.
The ellipsoid method has two key properties.
First, given a target accuracy $\epsilon$, it finds an $\epsilon$-approximate solution
in at most $O(d^2 \log(\epsilon^{-1}))$ subgradient queries.
Second, each iteration requires only $O(d^2)$ arithmetic operations to determine the next query.

We ask whether there is an algorithm for geodesically convex (g-convex) optimization on a Riemannian manifold which has linear convergence and each query can be computed efficiently. 
This problem has been stated informally before by~\citet[Sec.~2.2]{allenzhuoperatorsplitting} and~\citet[Sec.~1.1,1.3]{rusciano2019riemannian}.  Our contributions are to point out a partial solution for the case of constant curvature, to suggest possible avenues of attack, and to bring the problem to a wider audience.

Let us recall the relevant definitions.  {For references on g-convex optimization, see \citep{udriste1994convex} or~\citep[Ch.~11]{boumal2020intromanifolds}.}  
Throughout, $\calM$ denotes a complete $d$-dimensional Riemannian manifold with Riemannian metric $\langle \cdot, \cdot \rangle$, distance $\dist$, and exponential and logarithm maps $\exp$ and $\log$.
A subset $D$ of $\calM$ is \emph{g-convex} if for all $x, y \in D$
there is a unique minimizing geodesic segment $\gamma$ in $\calM$ connecting $x, y$, and $\gamma$ is contained in $D$ (this is sometimes called strongly g-convex~\citep[Sec.~11.3]{boumal2020intromanifolds}).
A function $f \colon D \rightarrow \reals$ is \emph{g-convex} in $D$ if
$f \circ \gamma$ is convex for all geodesic segments $\gamma$ contained in $D$.
A function $f$ is $M$-Lipschitz in $D$ if $|f(x) - f(y)| \leq M \dist(x, y)$ for all $x, y \in D$.
A tangent vector $g$ in the tangent space at $x$, denoted $\T_x \calM$, is a \emph{subgradient} of $f \colon D \to \reals$ at $x$ if $f(y) \geq f(x) + \langle g, \log_x(y)\rangle$ for all $y\in D$.  
A \emph{halfspace} on $\calM$ is a set $\{y \in \calM : \langle g, \log_x(y)\rangle \leq 0\}$ for some $(x, g)$ in the tangent bundle $\T\calM$.

We adopt the black-box model of optimization~\citep[Sec.~1.1]{nesterov2004introductory}.  A first-order algorithm can access the function $f$ to be minimized through oracle queries $x_k \in \calM$.
After query $x_k$, the algorithm is given $f(x_k)$ and a subgradient $g$ at $x_k$.
We can now state the open question:
\begin{question}\label{mainopenq}
Assume $\calM$ has sectional curvatures in the interval $[-K, K]$.
Assume $\bar B(\xorigin, r)$ is g-convex.  Let $f\colon \bar B(\xorigin, r) \to \reals$ be g-convex and $M$-Lipschitz in $\bar B(\xorigin, r)$.
Define $\zeta_{r \sqrt{K}} = \frac{r \sqrt{K}}{\tanh(r \sqrt{K})}$.
Is there a deterministic first-order algorithm $\calA$ with the following properties?
\begin{enumerate}
\item[(a)] For every $\epsilon \in (0,1)$, algorithm $\calA$ finds a point $x$ such that $f(x) - f^* \leq \epsilon \cdot M r$ in at most $O(\mathrm{poly}(\zeta_{r\sqrt{K}}, d) \log(\epsilon^{-1}))$ oracle queries, where $f^* = \min_{x \in \bar B(\xorigin, r)} f(x)$.
\item[(b)] Each iteration requires only a polynomial number of arithmetic operations, $\mathrm{poly}(\zeta_{r\sqrt{K}}, d)$, to determine the next query (in addition to the subgradient queries).
\end{enumerate}
\end{question}

Unlike in the Euclidean case, we permit the complexity of the method to depend on the curvature through $\zeta_{r \sqrt{K}} = \Theta(1 + r \sqrt{K})$ as such dependence is unavoidable~\citep{criscitiello2023curvature}.

A method satisfying properties (a) and (b) in Open Question~\ref{mainopenq} is of interest for several reasons.
First, the ellipsoid method is of fundamental theoretical importance in convex optimization.
As geodesic convexity is a generalization of convexity, it is natural to ask for such a generalization.  
Second, a method solving Open Question~\ref{mainopenq} applies to nonsmooth g-convex optimization problems, for example computing the geometric median such as for computational anatomy~\citep{Fletcher2009TheGM} or phylogenetics~\citep[Ch.~8]{bacak2014hadamard}.
Third, there are applications where the cost function is g-convex, but not strongly g-convex, and one seeks a linear convergence rate.
A notable example is \emph{operator scaling}~\citep{burgissernoncommutativeoptimization}, where it is also important that the method is deterministic.

The operator scaling problem encompasses several statistical problems including robust covariance estimation~\citep{wiesel2015gconvexity,sra2015conicgeometricoptimspd,franksmoitra2020} and estimation for matrix normal models ~\citep{tang2021integrated,franks2021neartyler}.
It also encompasses several questions in theoretical computer science, including a variant on polynomial identity testing (see references in~\citep{burgissernoncommutativeoptimization}).
For operator scaling, \citet{allenzhuoperatorsplitting} propose a linearly-convergent box-constrained Newton method.
This method does not solve Open Question~\ref{mainopenq} as it assumes additional properties about the objective (e.g., second-order robustness).


There are other methods for convex optimization which have similar properties as the ellipsoid method~\citep[pg.~156]{nesterov2004introductory}.
We focus on generalizing the ellipsoid method because it arguably has the simplest analysis among such methods, and similar obstacles appear in generalizing other methods.
In Section~\ref{constcurvature}, using the tools of geodesic maps introduced by~\citet{martinezrubio2021global}, we describe an ellipsoid-like method for spaces of constant curvature with query complexity $O(\zeta_{r \sqrt{K}} d^2 \log^2(\epsilon^{-1}))$ and per-query complexity $O(d^2)$.
In Section~\ref{obstacles}, we discuss possible approaches and obstacles to solving Open Question~\ref{mainopenq} for general Riemannian manifolds.

\smallskip \noindent 
\textbf{Related work} There are algorithms which satisfy properties (a) and (b) of Open Question~\ref{mainopenq} individually, but not concurrently.
Each iteration of the subgradient method~\citep{zhang2016complexitygeodesicallyconvex} can be computed efficiently.  However, without additional assumptions, this method does not have linear convergence.
The centerpoint method due to~\citet{rusciano2019riemannian} satisfies property (a).
However, each iteration requires computing a centerpoint, which we do not know how to compute efficiently.
\citet{lai2022riemannian} and~\cite{hirai2023interiorpoint} give Riemannian interior-point methods, which can converge quickly, but only under additional assumptions on the cost function.

\smallskip \noindent 
\textbf{Intermediate milestones} Open Question~\ref{mainopenq} asks for a deterministic algorithm, but one may initially permit the use of random queries.
For example, randomness can be helpful for computing the center of gravity in $\reals^d$~\citep[Sec.~3.3]{nemirovskibook}.
One may also focus on the Hadamard manifold of $n \times n$ positive definite matrices endowed with the affine-invariant metric~\citep[Sec.~11.7]{boumal2020intromanifolds}.  This is the most common manifold occurring in applications of g-convexity and possesses a well-studied structure~\citep{bridsonmetric,dolcetti2018differential}.  

\section{An ellipsoid-like method for spaces of constant curvature}\label{constcurvature}
Consider the setting described in Open Question~\ref{mainopenq}.
We first observe that it is essentially enough to solve Open Question~\ref{mainopenq} with radius $R = 1/\sqrt{K}.$
Let $D = \bar B(\xorigin, r)$ and assume $r > R$.  {(The case $r\leq R$ is easier to handle: we only need to pull $f$ back by a geodesic map once.)}
Consider the following algorithm.  
Initialize $x_0 = \xorigin$.
Given $x_k$, approximately solve the subproblem $\min_{\bar B(x_k, R)\cap D} f$ for a point $x_{k+1}$ such that $f(x_{k+1}) - \min_{\bar B(x_k, R) \cap D} f \leq \frac{\epsilon}{4} M R$.
Repeat the process.
Using only g-convexity of $f$ and $D$, one can show that after repeating this procedure $T = \ceil{\frac{2 r}{R} \log(\frac{2}{\epsilon})}$ times, we have $f(x_{T}) - f^* \leq \epsilon M r$ (see the proof of Theorem 7 in~\citep[App.~A]{martinezrubio2021global} for details).
Therefore, if we can solve each subproblem with $O(d^2 \log(\epsilon^{-1}))$ subgradient queries, then we can solve the original problem in $O(\zeta_{r\sqrt{K}} d^2 \log^2(\epsilon^{-1}))$ queries.

Now further assume $\calM$ is a $d$-dimensional hyperbolic space $\mathbb{H}^d$ (curvature equals $-1$ without loss of generality), and let us show how to solve each subproblem $\min_{\bar B(x_k, R)\cap D} f$.  
The analysis for a sphere is similar, so we omit it.
The key tool we need is geodesic maps.
A geodesic map with base point $x \in \calM$ is a diffeomorphism from $\mathbb{H}^d$ to the open unit ball $B(0,1)$ of $\reals^d$ which maps $x$ to the origin and which maps geodesics of $\mathbb{H}^d$ to straight lines of $\reals^d$ intersected with $B(0,1)$.
Geodesic maps of $\mathbb{H}^d$ have an explicit formula given by the Beltrami-Klein model, see~\citep[App.~C]{martinezrubio2021global}.
Our strategy to approximately solve $\min_{\bar B(x_k, R) \cap D} f$ is simple: we first pull $f$ back to a Euclidean space via a geodesic map $\phi_k \colon \mathbb{H}^d \to B(0,1)$ based at $x_k$, and then solve the resulting Euclidean problem $\min_{\bar B(0, \tilde R) \cap \phi_k^{-1}(D)} f \circ \phi_k^{-1}$ using the usual ellipsoid method.
We initialize that method with the ball $\bar B(0, \tilde R)$.  Here $\tilde R$ equals $R$ times an absolute constant, i.e., $\tilde R = \Theta(1/\sqrt{K})$.

Why does our strategy for solving the subproblem work?
First, geodesic maps have the key property 
that they map halfspaces in $\mathbb{H}^d$ to halfspaces of $\reals^d$ intersected with $B(0,1)$ (this fact can be verified by inspecting the explicit formula of a geodesic map).
In particular, given $y_k \in \bar B(0, \tilde R) \cap \phi_k^{-1}(D)$ (e.g., a query from the ellipsoid method), we can compute a Euclidean halfspace containing the sublevel set $\{y \in \phi_k^{-1}(D) \cap B(0,1) : f(\phi_k^{-1}(y)) \leq f(\phi_k^{-1}(y_k))\}$ by pulling back the corresponding Riemannian halfspace $\{x \in D : f(x) \leq f(\phi_k^{-1}(y_k))\}$ with $\phi_k$.\footnote{This is alternatively quantified with the notion of \emph{tilted convexity} introduced by~\citet[Lem.~3]{martinezrubio2021global}.}
Likewise if $y_k \not\in \bar B(0, \tilde R) \cap \phi_k^{-1}(D)$, one can compute a Euclidean separating halfspace for $\bar B(0, \tilde R) \cap \phi_k^{-1}(D)$ by pulling back a Riemannian separating halfspace for $\bar B(x_k, R) \cap D$.
Second, since we have chosen $R$ to be small, the metric distortion due to the geodesic map is small: in particular, $f \circ \phi_k^{-1}$ is still $O(M)$-Lipschitz~\citep[Lem.~2.c]{martinezrubio2021global}.
Lastly, note that the convexity of $f$ is not needed for the classical ellipsoid method to work: we just need that $f$ is Lipschitz and at each query point we are provided with a halfspace containing the minimizers of the problem (as guaranteed by our first observation).
One can verify that the Euclidean volume of $\bar B(0, \tilde R) \cap \phi_k^{-1}(D)$ is at least the volume of $B(0, c \tilde R)$ for some absolute constant $c \in (0,1)$.  This allows us to conclude the ellipsoid method finds an $\frac{\epsilon}{4} M R$-approximate solution to each subproblem in $O(d^2 \log(\epsilon^{-1}))$ queries.

Can we extend this technique beyond hyperbolic spaces?
Unfortunately, Beltrami's theorem states that the only Riemannian manifolds which admit geodesic maps to Euclidean space are those of constant curvature~\citep{beltramitheorem}.
One could replace a geodesic map with the exponential map, i.e., solve subproblems of the form $\min_{\T_x \calM} f \circ \exp_x$, 
and use comparison theorems~\citep{cheeger2008comparisontheorems} to analyze this method; 
but it is unclear how to make this approach work.
\section{Possible approaches and obstacles}\label{obstacles}
There is no convenient notion of ellipsoid on manifolds.
We consider two ways to handle this issue.



\smallskip \noindent 
\textbf{Ellipsoids in tangent spaces} 
One could maintain ellipsoids in the tangent spaces of $\calM$.
Consider the following algorithm.
Let $x_0 = \xorigin$ and $x^* \in \arg\min_{x \in \bar B(\xorigin, r)} f(x)$.
We initially know $\log_{x_0}(x^*)$ is contained in the ellipsoid $E_0 := \{v \in \T_{x_0} \calM : \|v\| \leq r\}$.
Querying at $x_0$ to get a subgradient $g_0$, we know $\log_{x_0}(x^*)$ is contained in the intersection $E_0 \cap \{v \in \T_{x_0} \calM : \langle v, g_0\rangle \leq 0\}$.  
Take an ellipsoid $\tilde E_0 \subset \T_{x_0} \calM$ with center $\tilde c_0$ which contains this intersection and has sufficiently small volume.  Define $x_1 = \exp_{x_0}(\tilde c_1)$.
To repeat this process, we need a method for \emph{transferring} the ellipsoid $\tilde E_0$ in $\T_{x_0} \calM$ to an ellipsoid $E_1$ in $\T_{x_1}\calM$ whose center is the origin, and such that $\exp_{x_1}(E_1) \supseteq \exp_{x_0}(\tilde E_0)$.  
\citet[Lem.~1,2]{kimyang} implicitly introduced a way of transferring \emph{balls} between tangent spaces while controlling the distortion of this transfer in order to obtain an accelerated method on manifolds.
It is not clear how to generalize their results to ellipsoids.

How should we keep track of the size of these ellipsoids?  We could track their Euclidean volumes.  However, even if the Euclidean volume of an ellipsoid $E \subset \T_x \calM$ is small, 
the Riemannian volume of $\exp_{x}(E)$ is not necessarily small, unless the ellipsoid is \emph{bounded}.
A bound on the Riemannian volume is important because it provides a bound on the optimality gap~\citep{rusciano2019riemannian}.

\smallskip \noindent 
\textbf{Halfspace- and volume-preserving maps}
A key ingredient in the analysis of the ellipsoid method is that there is a nontrivial family of \emph{halfspace-preserving} and \emph{volume-preserving} diffeomorphisms from $\reals^d$ to $\reals^d$: the invertible affine maps $x \mapsto A x + b$ with $\det(A) = 1$.
The explicit formula for the minimum-volume ellipsoid containing the intersection of an initial ellipsoid and a halfspace is derived by using a {volume-preserving} affine map to transform the ellipsoid to a ball (and the halfspace to another halfspace), and then computing the minimum-volume ellipsoid containing the intersection of a halfspace and a ball (a much easier calculation).
Let $\mathcal{G}$ be the group of halfspace-preserving and volume-preserving diffeomorphisms from $\calM$ to $\calM$.
Define the set of ``ellipsoids'' $\calE$ on $\calM$ as the result of applying all elements of $\mathcal{G}$ to every geodesic ball.
If $\calM$ is a Euclidean space, $\calE$ contains all ellipsoids in the usual sense.
Unfortunately, even for highly symmetric spaces like hyperbolic space, it seems $\mathcal{G}$ is the set of isometries of $\calM$~\citep{MOpost1,MOpost2}.
In this case, $\calE$ is the set of geodesic balls, which is not a sufficiently flexible family of subsets to be useful.

%



\bibliography{../bibtex/boumal}

\end{document}